\documentclass[12pt,twoside]{article}
\usepackage{amsmath}
\usepackage{hyperref}
\usepackage{amsfonts}
\usepackage{amssymb}

\usepackage{graphics}
\usepackage[pdftex]{graphicx}

\def\be{\begin{displaymath}}
\def\ee{\end{displaymath}}
\def\bee{\begin{equation}}
\def\eee{\end{equation}}
\def\Mas{Ma$\acute{\rm s}$lanka  }
\usepackage{graphicx}

\def\be{\begin{displaymath}}
\def\ee{\end{displaymath}}
\def\bee{\begin{equation}}
\def\eee{\end{equation}}
\def\Mas{Ma$\acute{\rm s}$lanka  }

\begin{document}
\title{On the Riesz and  B{\'a}ez-Duarte \\
criteria for the Riemann Hypothesis}
\author{
Jerzy Cis{\l}o, Marek Wolf \\
Institute of Theoretical Physics  \\
University of Wroc{\l}aw  \\
pl. M.Borna 9, 50-205 Wroc{\l}aw, Poland  \\
cislo@ift.uni.wroc.pl, mwolf@ift.uni.wroc.pl
}
\maketitle

\begin{abstract}
We investigate the relation between the Riesz and  the B{\'a}ez-Duarte criterion
for the Riemann Hypothesis. In particular we present the relation between the
function $R(x)$ appearing in the Riesz  criterion and the sequence $c_k$ appearing
in the B{\'a}ez-Duarte formulation. It is shown that $R(x)$ can be expressed by $c_k$,
and, vice versa, the sequence $c_k$ can be obtained from the
values of $R(x)$ at integer arguments. Also, we give some relations involving $c_k$ and $R(x)$,
and  value of the alternating sum of $c_k$.
\end{abstract}

\vspace{1cm}

\hfill{\it Dedicated to Prof. Luis B{\'a}ez-Duarte}

\hfill{\it on the occasion of his 70th birthday}\\

\bigskip

\section{Introduction}

The Riemann Hypothesis (RH) states that the nontrivial zeros of the function
\bee
\zeta(s)=\frac{1}{1-2^{1-s}}\sum_{n=1}^\infty \frac{(-1)^{n-1}}{n^s},
\label{zeta}
\eee
where $\Re(s)>0$ and $s\neq 1$ have the real part equal $\Re (s)=\frac{1}{2}$. Although Riemann
did not request it, today it is  often demanded additionally  that zeros on the critical line should
be  simple. The function $\zeta(s)$
defined by (\ref{zeta}) can be continued analytically to the whole
complex plane without $s=1$ where $\zeta(s)$ has the simple pole
\cite{Titchmarsh}.
There are probably over 100 statements equivalent to RH, see eg. \cite{Titchmarsh},
\cite{aimath}, \cite{Watkins}.
At the beginning of the 20th century M. Riesz \cite{Riesz} considered the function
\bee
R(x) = x \sum_{k=0}^\infty \frac{(-1)^{k}x^{k}}{k!\zeta(2k+2)}.
\label{Riesz}
\eee
We present the plot of $R(x)$ in the Fig. 1.

In \cite{Riesz} Riesz stated the following condition for the Riemann Hypothesis \\

{\bf Riesz Criterion:}
\bee
RH \Leftrightarrow   R(x) = \mathcal{O}\left( x^{1/4+\epsilon}\right) ~~~~~~~~~{\rm for~ each~~}  \epsilon>0.
\label{Riesz criterion}
\eee

A few years ago L. B{\'a}ez-Duarte \cite{Luis2}, \cite{Luis3}  considered the sequence
of numbers $c_k$ defined as the forward differences of $1/\zeta(2j+2)$:
\bee
c_k=\sum_{j=0}^k {(-1)^j \binom{k}{j}\frac{1}{\zeta(2j+2)}}.
\label{ckmain}
\eee
The plot of $c_k$ is shown on the Fig.2.

B{\'a}ez-Duarte proved\\

{\bf B{\'a}ez-Duarte  Criterion:}
\bee
RH  \Leftrightarrow c_k={\mathcal{O}}(k^{-\frac{3}{4}+\epsilon})~~~~~~~~~{\rm for~ each~~}  \epsilon>0.
\label{criterion}
\eee \\
Also, Baez-Duarte proved in \cite{Luis3} that it is not possible to replace $\frac{3}{4}$ by
larger exponent, and that $\epsilon = 0$ implies that the zeros
of $\zeta(s)$ are simple. Next in \cite{Luis4} B{\'a}ez-Duarte has considered
replacing "continuous" criteria with "sequential" criteria in more general setting.

Although  the title of the Baez-Duarte paper  \cite{Luis3} was {\it A sequential
Riesz-like criterion for the Riemann Hypothesis} he did not pursue further relation
between $c_k$ and $R(x)$ to prove his criterion  (he has used the Mellin transform).

In this paper we will present direct proof of the equivalence of the Riesz Criterion
and B{\'a}ez-Duarte Criterion. Besides, we will write some properties of $R(x)$ and of $c_k$,
of two-parameter generalizations of $R(x)$ and of $c_k$ introduced in
\cite{Paris}, \cite{Merlini},\cite{Aslam} and \cite{Coffey}. We calculate also the alternating
sum of $c_k$ and state the conjecture about the special sum of the M{\"o}bius function.

\section{Proof of equivalence of the Riesz and B{\'a}ez-Duarte Criteria}
In this section we will show that for large arguments function $R(x)/x$ and the sequence $c_k$
behave in a similar way.

Let $E$ denote the shift operator
\begin{equation}
Ef(n)=f(n+1).
\end{equation}
With this notation we can rewrite (\ref{ckmain}) as
\bee
c_k = (1-E)^k f(0), {~~\rm where~~} f(j)= \frac{1}{\zeta(2j+2)}.
\label{sequence}
\eee
We have the formal identity
\begin{equation}\label{form}
e^{x(1-E)}=e^xe^{-xE}.
\end{equation}
Consequently we have the identity
\begin{equation}
\sum_{k=0}^{\infty} \frac{x^k}{k!} \sum_{j=0}^k \binom{k}{j}  (-1)^j f(j)
=e^x \sum_{k=0}^{\infty} \frac{ (-x)^k f(k) }{k!}.
\end{equation}
After substitution $f(k)=1/\zeta(2k+2)$ we get
\begin{equation}
\sum_{k=0}^{\infty} \frac{ x^k }{k!} c_k = \frac{e^x}{x}R(x).
\label{Polya}
\end{equation}
We may also observe the above relation while comparing
discrete and continuous physical models of diffusion.
For both models we expect similar properties
of the solutions.

Relation (\ref{Polya}) appears also in the Exercises 67-71 in part IV of
the book Polya and Szeg\"{o} \cite{Polya}.

The formulae in the further part of this paper will involve the M\"{o}bius function
\begin{equation}\label{Mobius}
\mu(n) =
\left\{
\begin{array}{ll}
1 & \mbox {if $ n =1 $} \\
0 & \mbox {if $n$ is divisible by a square of a prime} \\
(-1)^k & \mbox{if $n$ is a product of $k$ different primes}
\end{array}
\right.
\end{equation}
 Using the formula
\begin{equation} \label{ odwr }
\frac{1}{\zeta(s)}=\sum_{n=1}^{\infty}\frac{ \mu (n) }{n^s }
\end{equation}
we can rewrite $R(x)$ and $c_k$ in the suitable for us form
\begin{equation} \label{ R }
R(x) = x\sum_{k=0}^{\infty}\frac{(-x)^k}{k! \zeta(2k+2)} =
x\sum_{n=1}^{\infty} \frac{\mu(n)}{n^2} \exp(-x/n^2),
\end{equation}
\begin{equation}\label{ck}
c_k=\sum_{j=0}^k  \binom{ k }{ j} \frac{(-1)^j}{\zeta(2j+2)}=
\sum_{n=1}^{\infty} \frac{\mu(n)}{n^2}\left( 1-\frac{1}{n^2} \right)^k.
\end{equation}
We will also consider two  two-parameter generalizations introduced in \cite{Paris},
\cite{Aslam}, \cite{Merlini}:
\begin{equation}
R_{a b} (x) = x \sum_{k=0}^{\infty}\frac{(-x)^k}{k! \zeta(ak+b)} =
\sum_{n=1}^{\infty} \frac{\mu(n)}{n^b} \exp(-x/n^a),
\end{equation}
\begin{equation}\label{ckab}
c_{a b}(k)=\sum_{j=0}^k  \binom{ k }{ j} \frac{(-1)^j}{\zeta(aj+b)}=
\sum_{n=1}^{\infty} \frac{\mu(n)}{n^b}\left( 1-\frac{1}{n^a} \right)^k.
\end{equation}
The original Riesz function $R(x)$ as well as the B{\'a}ez-Duarte sequence $c_k$ correspond to
the choice of parameters $a=b=2$. The generalization of the original Riesz criterion to
the family $R_{ab}(x)$ was given by A. Chaudhry \cite{Aslam}
\bee
RH \Leftrightarrow R_{ab}(x)=\mathcal{O} \left(x^{1-\frac{1}{a}(b-\frac{1}{2})+\epsilon}
\right)~~~~{\rm for~~ each} ~\epsilon>0.
\eee
For $a=2, ~b=1$ it reproduces the Hardy--Littlewood criterion \cite{Hardy} for RH.

We start from following simple lemma: \\

{\bf Lemma 1.} If  the function $f$ is nondecreasing for $1 \leq x \leq x_0$ and
nonincreasing for $x\geq x_0$ then
\begin{equation}
\sum_{n=1}^{\infty} f(n) \leq \int_1^{\infty} f(x) \;dx \; + f(x_0).
\end{equation} \\

Proof. If $f(1)\leq f(2) \leq \cdots \leq f(k)$ and $f(k) \geq f(k+1) \geq \cdots$, then
\bee
f(1)+f(2)+\cdots + f(k-1) \leq \int_1^k f(x)dx,
\eee
$$
f(k+1)+f(k+2)+f(k+3) \cdots \leq \int_k^{\infty} f(x)dx,
$$
$$
f(k) \leq f(x_0). ~~~~\Box
$$ \\

{\bf Corollary 1.} For $b>1$, $a >0 $, and $x>0$,  we have
\begin{equation}
\sum_{n=1}^{\infty} \frac{1}{n^b}\;\exp(-x/n^a) \leq
J_{ab} x^{(1-b)/a}   + \left( \frac{b}{ea} \right)^{b/a} x^{-b/a},
\end{equation}
where
\bee
J_{ab}=\int_0^{\infty}\frac{1}{t^b}\;\exp(-1/t^a)\; dt =  \frac{1}{a}\Gamma\left(\frac{b-1}{a}\right)
\eee
In particular we have
\bee
J_{2,2}=(1/2)\sqrt{\pi}, \quad J_{2,4}=(1/4)\sqrt{\pi}, \quad
J_{2,6}=(3/8)\sqrt{\pi}, \quad J_{2,8}=(15/16)\sqrt{\pi}.
\eee \\

{\bf Corollary 2.} We have
\begin{equation}
R_{ab}(x) =\mathcal{O}( x^{(1+a-b)/a} ),
\end{equation}
\bee \label{ckk}
c_{ab}(k)=\mathcal{O}( k^{(1-b)/a}.)
\eee
In particular, for $a=b=2$  we have:
\bee
\quad |R(x)| \leq  (1/2)\sqrt{\pi} x^{1/2} + 1/e.
\eee
The relation (\ref{ckk}) follows from the next lemma. \\

{\bf Lemma 2.} We have
\begin{equation}\label{Lemat}
\frac{R_{ab}(k)}{k} = c_{ab}(k)  + \mathcal{O}(k^{(1-a-b)/a}).
\end{equation} \\

Proof. For $x \in \langle 0, 1\rangle$, we have two inequalities:
\begin{equation}
\exp(-x) \geq 1-x, \quad \exp(x) \geq 1+ x +x^2/2.
\end{equation}
The first iequality implies
\begin{equation}\label{ineq1}
0\leq \exp(-kx) -(1-x)^k.
\end{equation}
The second inequality and Bernouli's inequality imply
\begin{equation}
(1-x)^k\exp(kx)\geq (1-x^2/2-x^3/2)^k \geq 1-kx^2/2-kx^3/2.
\end{equation}
After  some manipulations  we get
\begin{equation}\label{ineq2}
\exp(-kx)-(1-x)^k \leq (k/2)(x^2 + x^3) \exp(-kx).
\end{equation}
The inequalities (\ref{ineq1}) and  (\ref{ineq2}) give us estimation
\begin{equation}\label{main}
\left| \frac{R_{ab}(k)}{k} -c_{ab}(k) \right| \leq \sum_{n=1}^{\infty}
\frac{ \exp(-k/n^a) -  (1-1/n^a)^k }{n^b}
\end{equation}
$$
\leq \frac{k}{2}\sum_{n=1}^{\infty} \frac{1}{n^{2a+b}}\exp(-k/n^a)
 +   \frac{k}{2}\sum_{n=1}^{\infty} \frac{1}{n^{3a+b}}\exp(-k/n^a).
$$
Here we used the triangle inequality, the inequality $|\mu(n)| \leq 1$,
and the substitution $x=1/n^a$. Now the thesis of Lemma 2 follows from Corollary 1. $~~\Box$

The substitution $a=b=2$ in inequality (\ref{main}), and Corollary 1 give \\

{\bf Lemma 3.}
\begin{equation}
\left| \frac{R(k)}{k} -c_k \right| \leq \frac{3\sqrt{\pi}}{16} k^{-3/2} + \mathcal{O}(k^{-2}).
\label{bound}
\end{equation} \\

More explicitly we have
\begin{equation}
\left|\frac{R(k)}{k} -c_k \right| \leq \frac{3}{16}\sqrt{\pi} k^{-3/2} + \frac{27}{2} e^{-3} k^{-2}
+\frac{15}{32} \sqrt{\pi} k^{-5/2}  + 128 e^{-4} k^{-3}.
\label{R_k_c_k}
\end{equation}
Actually for $k>16$ we have $|R(k)/k-c_k| \leq (3/16)\sqrt{\pi} k^{-3/2}$.
Another proof of (\ref{R_k_c_k}) can be found in \cite{CW}, see also
\cite{wrzesien}.
The fact that approximately
$c_k \approx R(k)/k$
was observed previously by S. Beltraminelli and  D. Merlini \cite{Merlini}.

The Fig.3 shows depends on $k$ of $|R(k)/k-c_k|$ obtained
on the computer. Here the fit was obtained
by the least square method  from the data  with $k>10000$ to avoid transient regime
and it is given by the equation $y=0.01175x^{-1.527}$. \\

{\bf Lemma 4. } There is a real number $A$ such that for $0<x<y$
\begin{equation}\label{Rem}
\left| \frac{R(x)}{x} - \frac{R(y)}{y} \right|   \leq A (y-x) x^{-3/2}.
\end{equation} \\

Proof. We have
\begin{equation}
\left| \frac{R(x)}{x} - \frac{R(y)}{y} \right|
\leq \sum_{n=1}^{\infty} \frac{ \exp(-x/n^2) - \exp(-y/n^2)}{n^2}.
\end{equation}
From Mean-Value Theorem we conclude that there exists $z \in ( x, y )$ such that
\begin{equation}
\exp(-x/n^2) - \exp(-y/n^2) = \frac{y-x}{n^2} \exp(-z/n^2) < \frac{y-x}{n^2} \exp(-x/n^2).
\end{equation}
Finally it follows from Corollary 1 that
\begin{equation}
\left| \frac{R(x)}{x} - \frac{R(y)}{y} \right|    \leq (y-x) \left(\frac{\sqrt{\pi}}{4} x^{-3/2} +
\frac{4}{e^2}x^{-2} \right). ~~~\Box
\end{equation} \\

In paper \cite{Wolf}, the following equivalence had already been anticipated:  \\

{\bf Theorem 1.} For any real number $ \delta > -3/2 $ we have
\begin{equation} \label{Th}
R(x) = O(x^{\delta+1}) \Leftrightarrow  c_k = O(k^\delta).
\end{equation} \\

Proof. $\Rightarrow$ For integer $x$ (\ref{Th})  follows immediately from Lemma 3.
$\Leftarrow$ For non-integer $x$, we take Lemma 4 putting  $y=\lfloor x\rfloor+1$
and use (\ref{bound}). $~~\Box$ \\

{\bf Remark.} Putting $\delta=-3/4 + \epsilon $, we see that the Riesz criterion
is equivalent to the B{\'a}ez-Duarte criterion.

\section{The values of $c_k$ for large $k$}

For large negative $x$ function $R(x)$ tends to $xe^{-x}$.
For positive $x$, the behaviour of $R(x)$  is much more difficult to reveal
because the series (\ref{Riesz}) is very slowly convergent. Having applied
Kummer's acceleration convergence method, we  get
\bee
R(x) = x \left(\frac{6}{\pi^2} + \sum_{n=1}^\infty
\frac{\mu(n)}{n^2}\left(e^{-\frac{x}{n^2}} - 1\right)\right).
\label{Riesz3}
\eee
Using this formula we were able to produce the plot of $R(x)$ for $x$ up to $10^7$, see
Fig.1.  The first nontrivial zero of $R(x)$ is $x_0=1.156711643750816\ldots$. It is
a reflection of the fact, that $c_0>0$ while $c_1<0$.  Riesz in \cite{Riesz} has noticed the
existence of at least one  positive real zero of $R(x)$, while in \cite{Luis4} B{\'a}ez-Duarte
has proved existence of  infintely many zeros of $R(x)$.
The envelopes on the Fig.2  are given by the equations
\bee
y(x)=\pm A x^{\frac{1}{4}},~~~~ A=0.777506\ldots \times 10^{-5}.
\eee

It is very time consuming to calculate values of the sequence $c_k$ directly
from the definition (\ref{ckmain}), see \cite{Maslanka3}, \cite{Wolf}. The point is that for
large $j$, $\zeta(2j)$ is practically 1, and to distinguish it from 1 high precision
calculations are needed. The experience shows that to calculate
$c_k$ from (\ref{ckmain})  roughly $k\log_{10}(2)$ digits of accuracy is needed \cite{Wolf}.
However in \cite{Luis3} B{\'a}ez-Duarte gave the explicit formula for $c_k$ valid for large $k$:
\bee
c_{k-1}=\frac{1}{2k}\sum_{\rho} \frac{k^{\frac{\rho}{2}}\Gamma(1-\frac{\rho}{2})}{\zeta'(\rho)}
 + o(1/k),
\label{explicite}
\eee
where the sum runs over nontrivial zeros $\rho$ of $\zeta(s)$: $\zeta(\rho)=0$
and $\Im (\rho) \neq 0$.  \Mas in \cite{Maslanka3}
gives the similar formula which contains the term hidden in o(1/k) in (\ref{explicite}).
Let us introduce the notation
\bee
\frac{\Gamma(1-\frac{\rho_i}{2})}{\zeta'(\rho_i)} =  A_i+ iB_i.
\eee
Assuming  that $\rho_i=\frac{1}{2}+i\gamma_i$, it can be shown that
$A_i$ and  $B_i$ very  quickly decrease to zero \cite{Paris}, \cite{Wolf}:
\bee
\left|\frac{\Gamma(1-\frac{\rho_i}{2})}{\zeta'(\rho_i)}\right| \sim e^{-\pi\gamma_i/4}.
\label{decrease}
\eee
Finally, for large $k$, we obtain :
\bee
c_{k-1}=\frac{1}{k^{\frac{3}{4}}}\sum_{i=1}^\infty \left\{A_i \cos\left(\frac{\gamma_i \log(k)}{2}\right) -
B_i \sin\left(\frac{\gamma_i \log(k)}{2}\right)\right\}.
\label{rownanie}
\eee
The above formula explains oscillations on the plots of $c_k$ published in \cite{Luis3}
and \cite{Maslanka3}, see Fig.2. Because these curves are perfect cosine-like graphs
on the plots versus $\log(k)$ it means that in fact in the above formula (\ref{rownanie})
it suffices to maintain only  the first zero $\gamma_1=14.134725\ldots$, $A_1=2.0291739\ldots\times10^{-5},~~
B_1=-3.315924\ldots\times 10^{-5}$ and skip
all remaining  terms in the sum. It is justified by the very fast decrease of $A_i$ and $B_i$ following from
(\ref{decrease}).

\section{ The sums of $c_k$}

Let us perform  the formal calculation
\begin{eqnarray}
\sum_{k=0}^{\infty}t^k(1-E)^k & = &\frac{1}{1-t(1-E)}\\
= \frac{1}{1-t}~~\frac{1}{1+\frac{t}{1-t}E}& = & \frac{1}{1-t}\sum_{k=0}^{\infty}\left(- \frac{t}{1-t}\right)^kE^k. \nonumber
\end{eqnarray}
Acting with both sides  on the function $j \rightarrow 1/\zeta(2j+2)$ we get
\bee \label{sum}
\sum_{k=0}^{\infty}c_k t^k = \frac{1}{1-t}\sum_{k=0}^{\infty}\left(- \frac{t}{1-t}\right)^k
\frac{1}{\zeta(2k+2)}.
\eee
Of course  instead of $1/\zeta(2j+1)$ we can take an arbitrary function.

The above calculation was formal and we need to know what is the domain of convergence.
We will get (\ref{sum}) in another way.
Let us consider following identity
\bee
\frac{1}{n^2}\sum_{k=0}^{\infty}\left( 1 -\frac{1}{n^2} \right)^k t^k
=\frac{1}{t+(1-t)n^2} =
\frac{1}{1-t}\sum_{k=0}^{\infty}\left(- \frac{t}{1-t}\right)^k\frac{1}{n^{2k+2}}.
\eee
The first sum is convergent for $-1 \leq t  \leq 1$ while the second one is convergent for  $-\infty< t <1/2$.
Thus the common domain of convergence is the interval $\langle -1, 1/2)$. Hence
\bee \label{series}
\sum_{n=1}^{\infty}\sum_{k=0}^{\infty} \frac{\mu(n)}{n^2}\left( 1 -\frac{1}{n^2} \right)^k t^k
=\frac{1}{1-t}\sum_{n=1}^{\infty}\sum_{k=0}^{\infty} \frac{\mu(n)}{n^{2k+2}}\left(- \frac{t}{1-t}\right)^k.
\eee
The sums (\ref{series}) are absolutely convergent and we can change the order of summation
obtaining (\ref{sum}).

Substituting $t=-1$ in the equation (\ref{series}), we get
\bee
 \sum_{k=0}^{\infty}(-1)^k c_k  = \sum_{k=1}^{\infty}\frac{1}{2^k} \frac{1}{\zeta(2k)}
  =0.782527985325384234576688\ldots.
\label{suma}
\eee
This number probably can not be expressed by other known constants, because the Simon Plouffe
inverter failed to find any relation \cite{Plouffe}. Applying  Abel's summation, we can write the r.h.s.
of (\ref{suma}) as:
\bee
\sum_{k=1}^{\infty}\frac{1}{2^k} \frac{1}{\zeta(2k)} =
1+\sum_{k=1}^{\infty}\left( 1- \frac{1}{2^k}\right)\left( \frac{1}{\zeta(2k)}-\frac{1}{\zeta(2k+2)}\right)
\eee
$$
= 1 + \int_2^\infty \left(1-\frac{1}{2^{\lfloor x/2 \rfloor }} \right) \frac{\zeta'(x)}{\zeta^2(x)} dx.
$$
More detailed considerations gives
\bee
\sum_{j=0}^{k-1} (-1)^j c_j = \sum_{k=1}^{\infty}\frac{1}{2^k} \frac{1}{\zeta(2k)} - \frac{(-1)^k}{2} c_k
+ \mathcal{O}(k^{-3/2}).
\eee

Now we turn to the sum $\sum_{i=0}^{\infty} c_i$. The partial sum can be expressed
in the following way:
\bee
S_{k-1}=\sum_{i=0}^{k-1} c_i  =
\sum_{n=1}^{\infty} \mu(n) \left(1 - \left(1-\frac{1}{n^2}\right) ^k \right)=
-\sum_{j=1}^k \binom{k}{j}\frac{(-1)^j}{\zeta(2j)}.
\label{suma_c_k}
\eee
Computer calculations  show that the partial sums
initially tend from above to -2, but for $k\approx 91000$
the partial sum crosses -2 and around $k\approx 100000$ the partial sum starts to increase.
These oscillations begins to repeat with growing amplitude around -2, see Fig. 4.
The value -2 was informally derived in \cite{CW}.

For large $k$ the oscillations are described  by the  integral of  (\ref{rownanie})
$$
k^{1/4}\sum_{i=1}^\infty \frac{1}{1/4+\gamma_i^2}
\left\{(A_i+2B_i \gamma_i) \cos\left(\frac{\gamma_i \log(k)}{2}\right) -
 (B_i-2A_i \gamma_i)\sin\left(\frac{\gamma_i \log(k)}{2}\right)\right\}
$$
\bee
 = \mathcal{O}\left(k^{\frac{1}{4}}\right).
\label{k14}
\eee
It is interesting that the amplitude  is  very small, e.g. at
$k\sim 10^8$ the amplitude  is of the order 0.001. By combining (\ref{suma_c_k}) and
(\ref{k14}) we get that
\bee
\sum_{n=1}^{\infty} \mu(n) \left(1 - \left(1-\frac{1}{n^2}\right)^k \right)
\eee
oscillates around -2 with the amplitude growing like $k^{1/4}$. We generalize the last
statement in the form of the following

{\bf Conjecture 1:} Let $b\geq a>0$. Then the sum
\bee
\sum_{n=1}^{\infty} \frac{\mu(n)}{n^{b-a}} \left(1 - \left(1-\frac{1}{n^a}\right) ^k \right)
\eee
oscillate around
\bee
1/\zeta(b-a)
\eee
with the amplitude given unconditionally by $k^{\frac{1-b+2a}{2a}}$ and with the amplitude
growing like $k^{\frac{a-b+1/2}{a}}$ under the assumption of the Riemann Hypothesis.

It seems to be mysterious that the sum $\sum_{i=0}^{\infty} c_i$ oscillates around -2, while
the alternating sum $\sum_{i=0}^{\infty}(-1)^i c_i$ gives probably transcendent number.\\

{\bf Acknowledgement} We would like to thank Prof. L. B{\'a}ez-Duarte, Prof. M. Coffey  and Prof. K. \Mas
for e-mail exchange.
To prepare  data for some figures we used  the free computer algebra system PARI/GP
\cite{Pari}. \\

%\end{document}

\newpage

\begin{figure}[pht]
\vspace{-2.5cm}
\begin{minipage}{15.8cm}
\begin{center}
\hspace{-3.5cm}
\includegraphics[width=12cm,angle=0, scale=1]{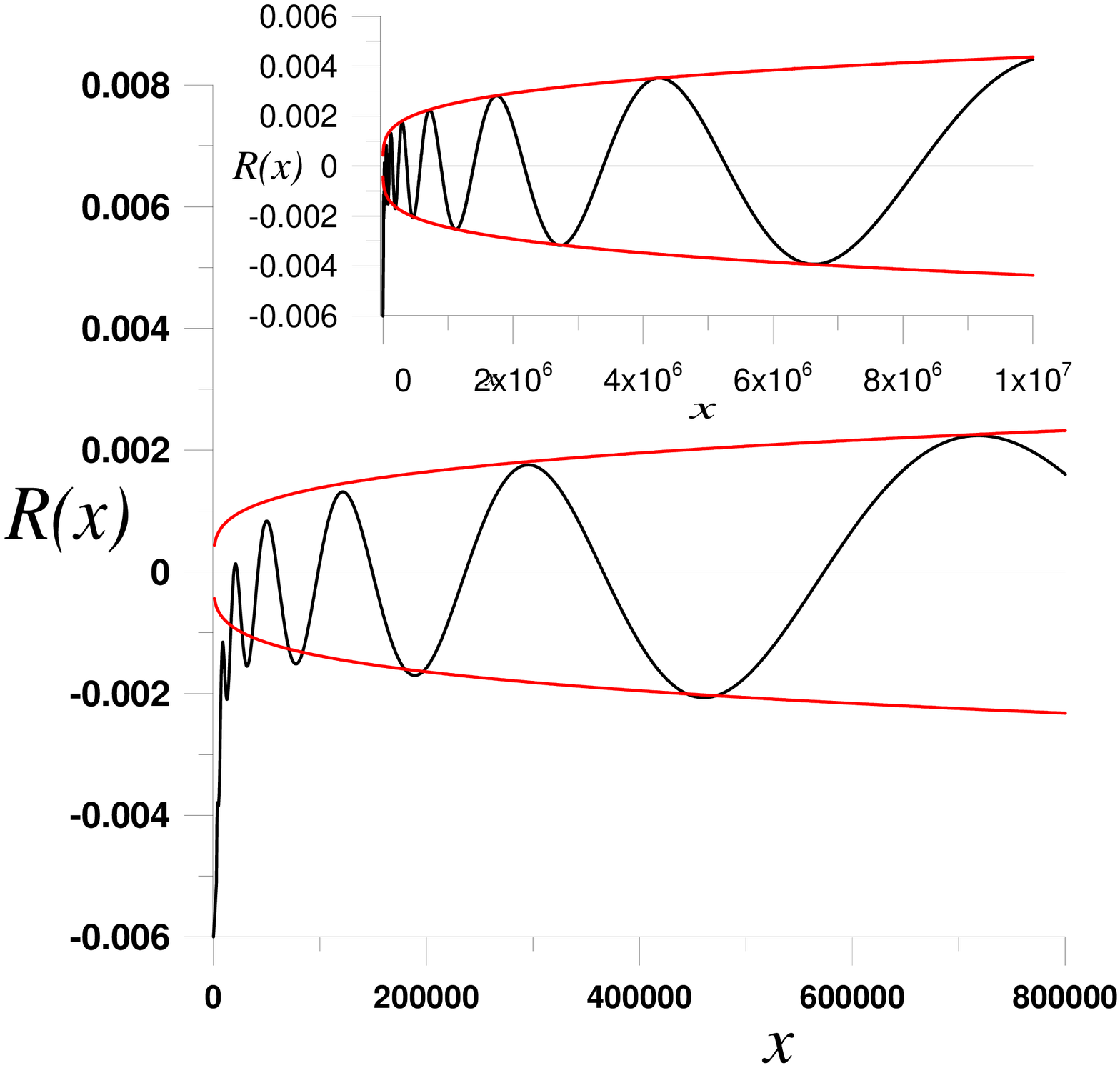} \\
\vspace{-3.5cm} Fig.1  The plot of $R(x)$ for $x\in(0, 800000)$ and for $x\in(0, 10^7)$ in
the inset. The part of $R(x)$ smaller than -0.006 is skipped.  \\
\end{center}
\end{minipage}

%\vspace{2.0cm}
\begin{minipage}{15.8cm}
\begin{center}
\hspace{-3.5cm}
\includegraphics[width=12cm,angle=0,angle=-90, scale=1]{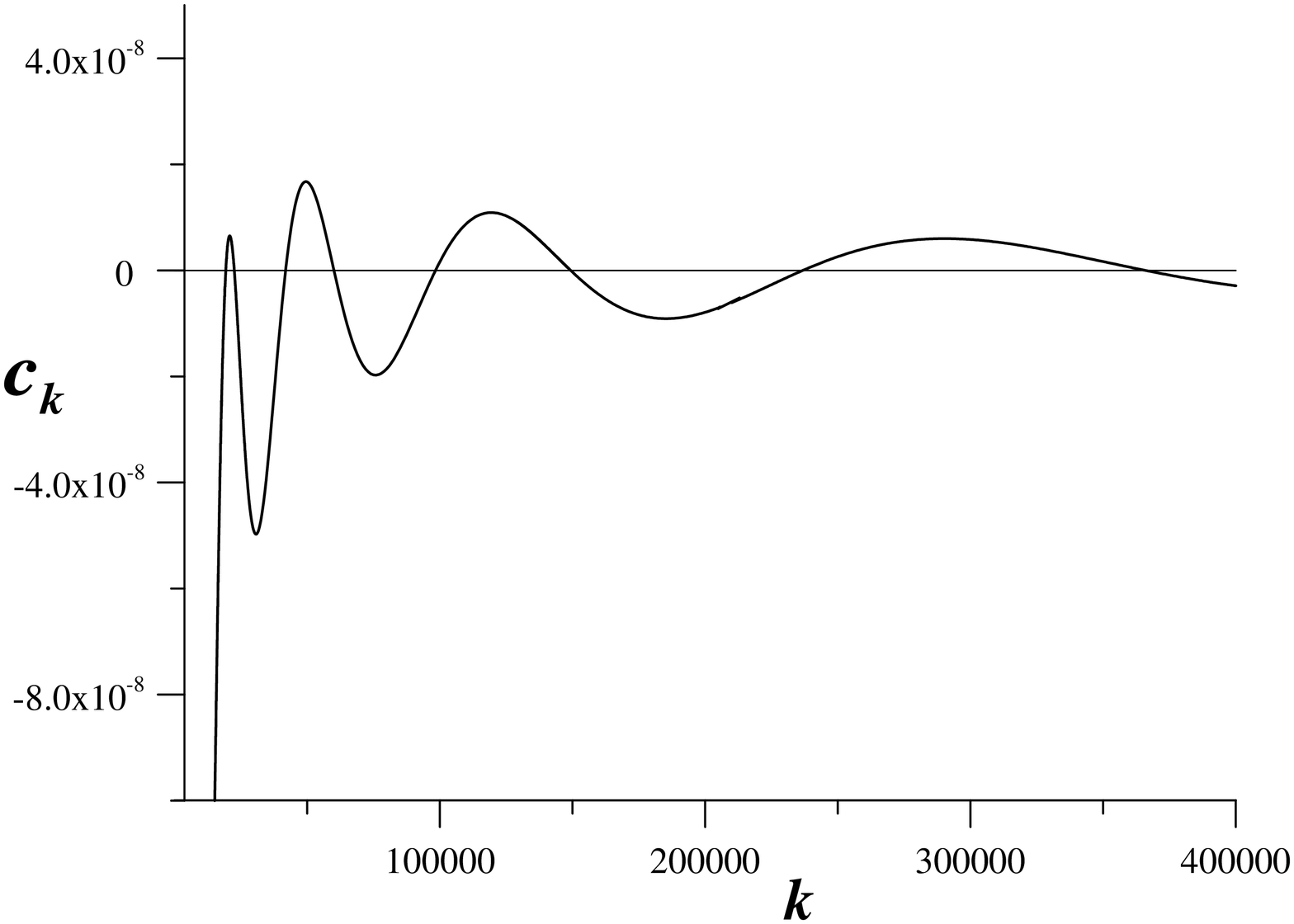} \\
Fig.2  The plot of $c_k$ for $k\in(1, 400000)$. \\
\end{center}
\end{minipage}
\end{figure}

\newpage

\begin{figure}
\begin{minipage}{12.8cm}
\begin{center}
\vspace{-2.0cm}
\includegraphics[width=10cm, scale=0.8]{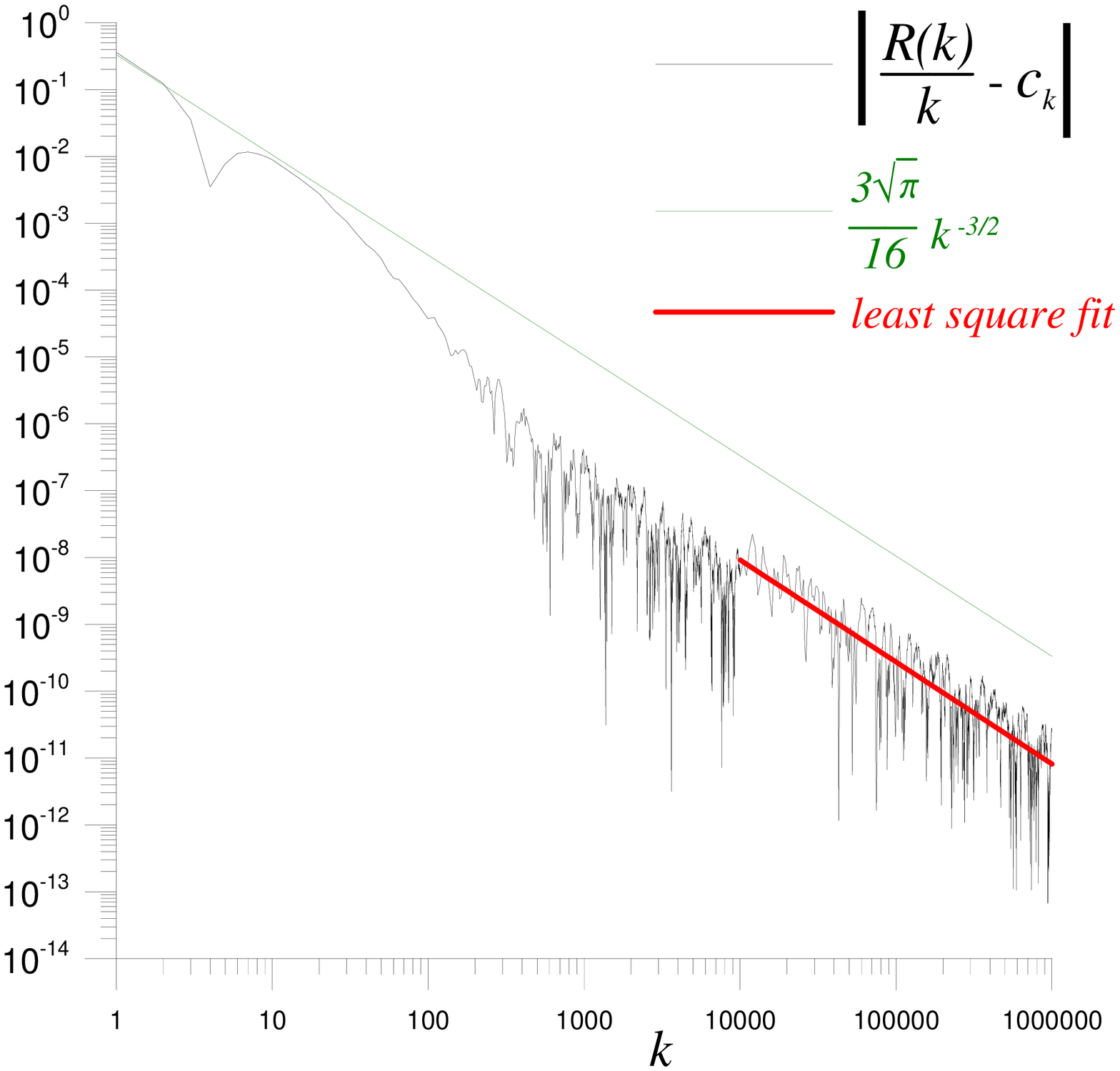}\\
\vspace{-2.5cm}
Fig. 3  The log-log plot of $|R(k)/k-c_k|$ for $k\in(0,10^6)$\\
\end{center}
\end{minipage}

\begin{minipage}{12.8cm}
\begin{center}
\vspace{-1.0cm}
\includegraphics[width=10cm, scale=1]{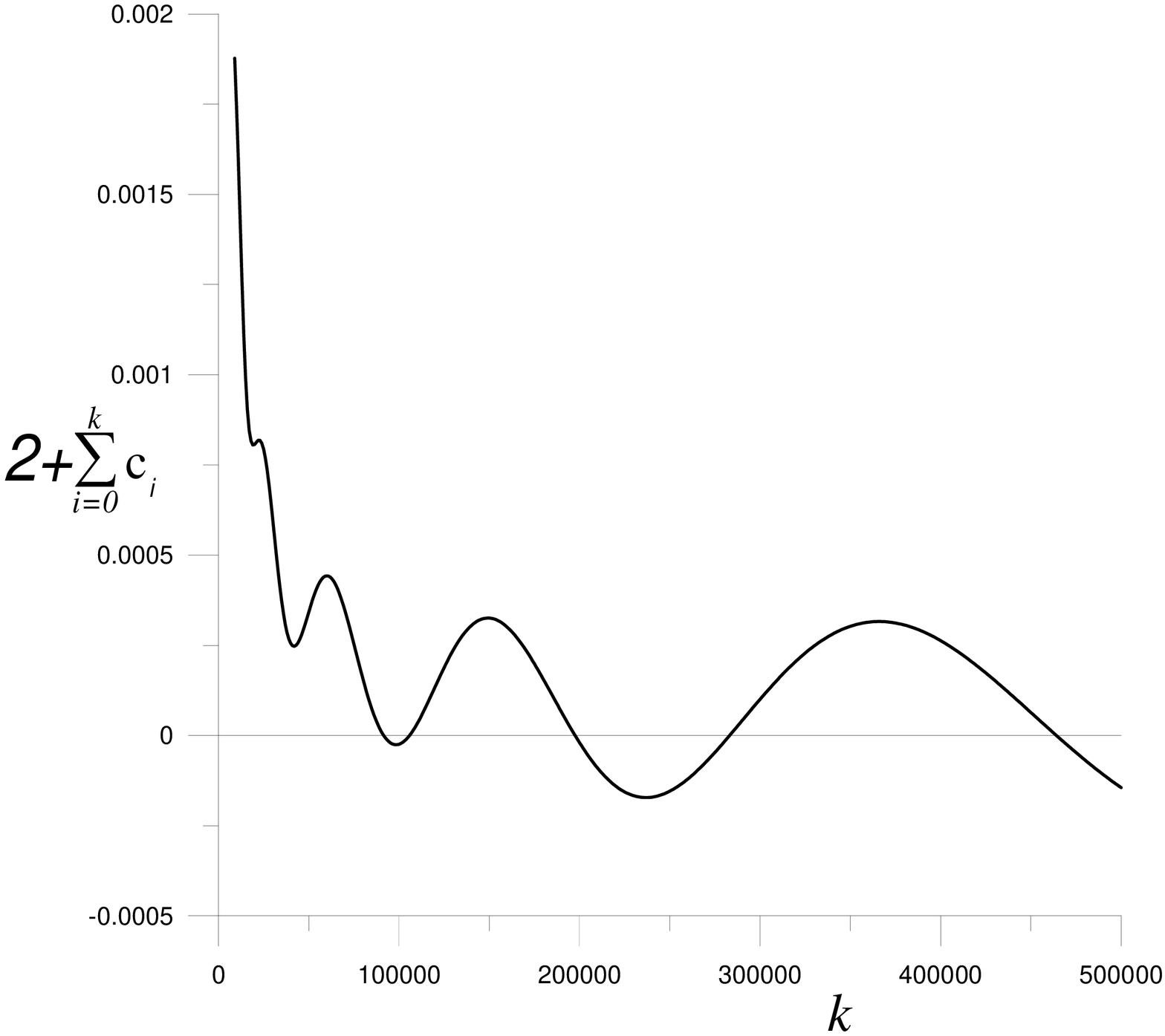}\\
\vspace{-2.5cm}
Fig. 4 The distance from -2 of the partial sums $\sum_{k=0}^n c_k$ for $n=1,\ldots 500000$.
\end{center}
\end{minipage}
\end{figure}

\end{document}